\documentclass[12pt]{article}
\usepackage{amssymb,amsmath,amsfonts,amsthm,amscd,latexsym,indentfirst,verbatim,wasysym}
\def\Var{\mathrm{Var}}
\def\Span{\mathrm{Span}}

\def\As{\mathrm{As}}

\begin{document}

\begin{flushright}
MSC (2010): 16W99 
\end{flushright}

\begin{center}
{\Large
Embedding of pre-Lie algebras into preassociative algebras}

V. Gubarev
\end{center}

\begin{abstract}
With the help of Rota---Baxter operators and the Gr\"{o}bner---Shirshov bases, 
we prove that every pre-Lie algebra injectively embeds into its universal 
enveloping preassociative algebra. 

\medskip
{\it Keywords}:
Rota---Baxter operator, Gr\"{o}bner---Shirshov basis, 
pre-Lie algebra, preassociative algebra (dendriform algebra).
\end{abstract}

\section{Introduction}

Pre-Lie algebras were introduced independently by E.B.~Vinberg \cite{Vinberg60} in 1960
and M.~Gerstenhaber \cite{Gerst63} in 1963. 
Pre-Lie algebras also known as left-symmetric algebras satisfy the identity 
$(x_1 x_2)x_3 - x_1(x_2 x_3) = (x_2 x_1) x_3 - x_2 (x_1 x_3)$.

In 2001, J.-L. Loday \cite{Dialg99} defined the dendriform (di)algebra
(preassociative algebra) as a vector space endowed 
with two bilinear operations $\succ,\prec$ satisfying 
$$
\begin{gathered}
(x_1\succ x_2+x_1\prec x_2)\succ x_3 = x_1\succ (x_2 \succ x_3), \
(x_1\succ x_2)\prec x_3=x_1\succ(x_2\prec x_3), \\
x_1\prec(x_2\succ x_3+x_2\prec x_3)=(x_1\prec x_2)\prec x_3.
\end{gathered}
$$

In 1995, J.-L. Loday also defined~\cite{Loday95} Zinbiel algebra
(precommutative algebra), on which the identity
$(x_1x_2 + x_2x_1)x_3 = x_1(x_2 x_3)$ holds. 
Every preassociative algebra with the identity $x\succ y = y\prec x$ 
is a precommutative algebra ($x_1x_2 = x_1\succ x_2$)
and under the product $x\cdot y = x\succ y - y\prec x$ is a pre-Lie algebra.

A linear operator $R$ defined on an algebra $A$ over a field $\Bbbk$
is called a Rota---Baxter operator (RB-operator) of a weight $\lambda\in\Bbbk$
if it satisfies the relation
\begin{equation}\label{RB}
R(x)R(y) = R( R(x)y + xR(y) + \lambda xy), \quad x,y\in A.
\end{equation}
An algebra $A$ with a Rota---Baxter operator is called a Rota---Baxter algebra 
(RB-algebra)~\cite{Baxter60,Guo2011}.

In 2000, M. Aguiar \cite{Aguiar00} stated that an associative 
(commutative) algebra with a given Rota---Baxter operator~$R$
of weight zero under the operations $a \succ b = R(a)b$, $a\prec b = aR(b)$ 
is a preassociative (precommutative) algebra. 
The analogue of this construction for the pair of pre-Lie algebras 
and Lie RB-algebras of weight zero was stated 
in 2000 by M.~Aguiar~\cite{Aguiar00} 
and by I.Z. Golubchik, V.V. Sokolov~\cite{GolubchikSokolov}.
In 2013 \cite{BBGN2012}, the construction of M.~Aguiar 
was generalized for the case of an arbitrary variety.

In 2008, the notion of universal enveloping RB-algebras
of preassociative algebras was introduced by K. Ebrahimi-Fard and L. Guo~\cite{FardGuo07}.
Moreover, it was proved in~\cite{FardGuo07} that the universal enveloping of 
a free preassociative algebra is free.

In 2010, with the help of Gr\"{o}bner---Shirshov bases \cite{BokutChen},
Yu. Chen and Q.~Mo stated that a preassociative algebra
over a field of characteristic zero injectively embeds into
its universal enveloping RB-algebra \cite{Chen11}.

In 2013 \cite{GubKol2013}, given a variety $\Var$, it was proved that every pre-$\Var$-algebra
injectively embeds into its universal enveloping $\Var$-RB-alge\-bra of weight zero. 

{\bf Problem 1.1}.
Prove that every pre-Lie algebra injectively embeds into
its universal enveloping preassociative  algebra.

{\bf Problem 1.2}~\cite{Kol2017}.
Construct the universal enveloping preassociative algebra for a given pre-Lie algebra. 

A solution of Problem~1.2 will cover Problem~1.1.  
The discussion of Problems~1,\,2 in the case of restricted pre-Lie algebras can be found in \cite{Dokas13}. 
The analogues of Problems~1,\,2 for Koszul-dual objects, dialgebras, were solved in \cite{LP93}.

The main goal of the current work is to solve Problem~1.1.
In September 2018, the proof that the pair of varieties of pre-Lie
and preassociative algebras is a PBW-pair (see~\cite{PBW}) appeared in~\cite{Dots},
so it covers the solution of Problem~1.1. 
Anyway, the methods used are quite different. 
V. Dotsenko and P. Tamaroff are applying in~\cite{Dots} the operad theory, 
while the author uses Gr\"{o}bner---Shirshov bases technique developed for associative Rota---Baxter algebras.

In~\cite{Gub2018}, the analogue of Problem~1.1 for postLie and postassociative algebras was solved,
see there the more detailed history and description of the problems.

\section{Embedding of prealgebras into RB-algebras}

{\bf Theorem 2.1}~\cite{Aguiar00,BBGN2012,GolubchikSokolov}.
Let $A$ be an RB-algebra of a variety $\Var$ and weight zero.
With respect to the operations
\begin{equation}\label{LodayToRB}
x\succ y = R(x)y,\quad x\prec y = xR(y)
\end{equation}
$A$ is a pre-$\Var$-algebra.

Denote the pre-$\Var$-algebra obtained in Theorem~2.1 as $A^{(R)}$. 

Given a pre-$\Var$-algebra $\langle C,\succ,\prec\rangle$, universal enveloping 
RB-$\Var$-algebra $U$ of $C$ is the universal algebra in the class of all RB-$\Var$-algebras
of weight zero such that there exists a~homomorphism from $C$ to $U^{(R)}$.

{\bf Theorem 2.2}~\cite{GubKol2013}.
Any pre-$\Var$-algebra could be embedded into its universal enveloping 
RB-algebra of the variety $\Var$ and weight zero.

Let us consider the idea of the proof of Theorem~2.2.
Suppose $\langle A,\succ,\prec\rangle$ is a pre-$\Var$-algebra.
Then the direct sum of two isomorphic copies of $A$,
the space $\hat A = A\oplus A'$, endowed with a binary operation
\begin{equation}\label{hat-construction}
a*b = a\succ b + a\prec b, \quad
a*b' = (a\succ b)', \quad
a'*b = (a\prec b)', \quad
a'*b' = 0,\quad a,b\in A,
\end{equation}
is proved to be an algebra of the variety $\Var$.
Moreover, the map $R(a') = a$, $R(a) = 0$ is 
an RB-operator of weight~zero on $\hat A$. 
The injective embedding of $A$ into $\hat A$ is given 
by $a\mapsto a'$, $a\in A$.

\section{Embedding of pre-Lie algebras into preassociative algebras}

Let $R\As\langle X\rangle$ denote the free associative algebra 
generated by a set $X$ with a linear map $R$ in the signature.
One can construct a linear basis of $R\As\langle X\rangle$ (see, e.g., \cite{Guo2013})
by induction. At first, all elements from $S(X)$, the free semigroup generated by~$X$, 
lie in the basis. At second, if we have basic elements $a_1,a_2,\ldots,a_k$, $k\geq1$, 
then the word $w_1R(a_1)w_2\ldots w_kR(a_k)w_{k+1}$ lies in the basis 
of~$R\As\langle X\rangle$.
Here $w_2,\ldots,w_k\in S(X)$ and $w_1,w_{k+1}\in S(X)\cup \emptyset$,
where $\emptyset$ denotes the empty word.
Let us denote the basis obtained as $RS(X)$.
Given a word~$u$ from $RS(X)$, the number of appeareances of the symbol~$R$
in~$u$ is denoted by $\deg_R(u)$, the $R$-degree of~$u$. 
We call an element from $RS(X)$ of the form $R(w)$ as $R$-letter.
By~$X_\infty$ we denote the union of $X$ and the set of all $R$-letters.
Given $u\in RS(X)$, define $\deg u$ (degree of $u$) as the length of $u$
in the alphabet~$X_\infty$.

Suppose that $X$ is a well-ordered set with respect to $<$.
Let us introduce by induction the deg-lex order on $S(X)$. 
At first, we compare two words $u$ and $v$ by the length:
$u < v$ if $|u|<|v|$. At second, when $|u| = |v|$, 
$u = x_i u'$, $v = x_j v'$, $x_i, x_j\in X$,
we have $u < v$ if either $x_i < x_j$ or $x_i = x_j$, $u'< v'$.
We compare two words $u$ and $v$ from $RS(X)$ 
by $R$-degree: $u<v$ if $\deg_R(u)<\deg_R(v)$. 
If $\deg_R(u) = \deg_R(v)$, we compare $u$ and $v$ in deg-lex order as words 
in the alphabet $X_\infty$. Here we define each $x$ from $X$ to 
be less than all $R$-letters and $R(a)<R(b)$ if and only if $a<b$.

Let $*$ be a symbol not containing in $X$.
By a $*$-bracketed word on $X$, we mean a basic word 
from $R\As\langle X\cup\{*\}\rangle$ with exactly one occurrence of $*$. 
The set of all $*$-bracketed words on $X$ is denoted by $RS^*(X)$.
For $q\in RS^*(X)$ and $u\in R\As\langle X\rangle$, we define
$q|_u$ as $ q|_{*\to u}$ to be the bracketed word obtained by 
replacing the letter $*$ in $q$ by $u$.

The order defined above is {\it monomial}, i.e., from $u < v$  
follows that $q|_u < q|_v$ for all $u,v\in RS(X)$ and $q\in RS^*(X)$.

Given $f\in R\As\langle X\rangle$, by $\bar{f}$ we mean the leading word in $f$.
We call $f$ monic if the coefficient of $\bar{f}$ in $f$ is~1. 

{\bf Definition 3.1}~\cite{Guo2013}.
Let $f,g\in R\As\langle X\rangle$.
If there exist $\mu,\nu,w\in RS(X)$ such that 
$w = \bar{f} \mu = \nu \bar{g}$ with $\deg w<\deg(\bar{f})+\deg(\bar{g})$, 
then we define $(f,g)_w$ as $f\mu - \nu g$ and call it 
the {\it composition of intersection} of $f$ and $g$ with respect to $(\mu,\nu)$.
If there exist $q\in RS^*(X)$ and $w\in RS(X)$ 
such that $w = \bar{f} = q|_{\bar{g}}$, then we define
$(f,g)^q_w$ as $f-q|_g$ and call it the {\it composition of inclusion} 
of $f$ and $g$ with respect to $q$.

{\bf Definition 3.2}~\cite{Guo2013}.
Let $S$ be a subset of monic elements from $R\As\langle X\rangle$ 
and $w\in RS(X)$. 

(1) For $u,v\in R\As\langle X\rangle$, we call $u$ and $v$ 
congruent modulo $(S, w)$ and denote this by $u \equiv v \mod (S, w)$ 
if $u - v = \sum c_i q_i|_{s_i}$ 
with $c_i \in \Bbbk$, $q_i\in RS^*(X)$, $s_i\in S$ and $q_i|_{\overline{s_i}} < w$.

(2) For $f,g\in R\As\langle X\rangle$ and suitable 
$w,\mu,\nu$ or $q$ that give a composition of intersection
$(f,g)_w$ or a composition of inclusion $(f,g)^q_w$, 
the composition is called trivial modulo $(S, w)$ if
$(f,g)_w$ or $(f,g)^q_w \equiv 0 \mod (S, w)$.

(3) The set $S\subset R\As\langle X\rangle$ is 
called a {\it Gr\"{o}bner---Shirshov basis} if, for all $f,g\in S$, 
all compositions of intersection $(f,g)_w$ 
and all compositions of inclusion $(f,g)^q_w$ are trivial modulo $(S, w)$.

\newpage

{\bf Theorem 3.3}~\cite{Guo2013}. 
Let $S$ be a set of monic elements in $R\As\langle X\rangle$, 
let $<$ be a monomial ordering on $RS(X)$ and let $Id(S)$ 
be the $R$-ideal of $R\As\langle X\rangle$ generated by $S$. 
If $S$ is a Gr\"{o}bner---Shirshov basis in $R\As\langle X\rangle$,
then $R\As\langle X\rangle = \Bbbk Irr(S)\oplus Id(S)$ 
where $Irr(S) = RS(X)\setminus \{q|_{\bar{s}}\mid q\in RS^*(X),s\in S\}$
and $Irr(S)$ is a linear basis of $R\As\langle X\rangle/Id(S)$.

Let $A$ be an associative algebra with an RB-operator $R$.
Then the algebra $A^{(-)}$ is a Lie RB-algebra
under the product $[x,y] = xy-yx$ and the same action of $R$. 
Thus, we can state the analogues of Problems~1 and~2 for 
the varieties of Lie and associative RB-algebras. 
(1) Whether a Lie RB-algebra can always be embedded into an associative RB-algebra?
(2) What is a linear basis of the universal enveloping associative RB-algebra 
for a given Lie RB-algebra? 

Let $\hat{L} = L\oplus L'$ be exactly the Lie algebra 
with the RB-operator $R$ of weight~0 constructed 
in the sketch of the proof of Theorem~2.2. 
Let us choose linear bases $x_\alpha$ and
$y_\alpha$, $\alpha\in\Lambda$, of $L$ and $L'$ 
respectively such that
\begin{equation}\label{R-Cond}
R(y_\alpha) = x_\alpha,\quad  
R(x_\alpha) = 0.
\end{equation}
Then the set
$X_\Lambda\cup Y_\Lambda = \{x_\alpha,\alpha\in\Lambda\}\cup \{y_\alpha,\alpha\in\Lambda\}$
is a linear basis of $\hat{L}$.
Suppose that $\Lambda$ is a well-ordered set. 
Extend the order on the set $X_\Lambda\cup Y_\Lambda$ in the following way:
$t_\alpha<t_\beta$, $t\in\{x,y\}$, if and only if $\alpha<\beta$;
$y_\alpha<x_\beta$ for all $\alpha\leq \beta$
and $x_\beta<y_\alpha$ if $\beta<\alpha$.

Given $x\in X_\Lambda$, $y\in Y_\Lambda$, 
denote $[[\ldots[[y,x],x]\ldots],x]\in \hat{L}^{p+1}$ by $[y,x^{(p)}]$. 

Consider the set $S$ of the following elements 
in $R\As\langle X_\Lambda\cup Y_\Lambda\rangle$:
\begin{equation}
\begin{gathered}\label{UnivRel}
x_{\alpha}x_\beta - x_\beta x_{\alpha} - [x_\alpha,x_\beta],\ 
y_{\alpha}y_\beta - y_\beta y_{\alpha} - [y_\alpha,y_\beta],\ \beta<\alpha,\\
x_{\alpha}y_\beta - y_\beta x_{\alpha} - [x_\alpha,y_\beta],\ \beta\leq \alpha,\quad
y_\beta x_{\alpha} - x_{\alpha}y_\beta - [y_\beta,x_\alpha],\ \alpha<\beta,
\end{gathered}
\end{equation}

\vspace{-0.5cm}

\begin{gather}
R(a)R(b) - R(R(a)b + aR(b)), \label{almostRB} \\
R(R(z_1)\vec{x}_{\alpha_1}R(z_2)\ldots R(z_s)\vec{x}_{\alpha_s}R(z_{s+1})), \label{XRel} 
\end{gather}

\vspace{-0.8cm}

\begin{multline}\label{LongRel}
R\big(R(z_1)\vec{x}_{\alpha_1} R(z_2)\ldots R(z_s)\vec{x}_{\alpha_s}y_\beta x_\beta^k R(z_{s+1})\big) 
 - \frac{1}{k+1}\bigg(R(z_1)\ldots R(z_s)\vec{x}_{\alpha_s} x_\beta^{k+1} R(z_{s+1}) \\
 + \sum\limits_{i=2}^{k+1}(-1)^i\binom{k+1}{i}R\big(R(z_1)\vec{x}_{\alpha_1}R(z_2)\ldots R(z_s)\vec{x}_{\alpha_s}\big[y_\beta,x_\beta^{(i-1)}\big]x_\beta^{k+1-i} R(z_{s+1})\big) \\
 - \sum\limits_{j=1}^{s+1}R\big(R(z_1)\vec{x}_{\alpha_1}R(z_2)\ldots R(z_{j-1})\vec{x}_{\alpha_{j-1}}z_j\vec{x}_{\alpha_j}R(z_{j+1})\ldots R(z_s)\vec{x}_{\alpha_s}x_\beta^{k+1}R(z_{s+1})\big) \\
 - \sum\limits_{j=1}^{s}\sum\limits_{t=1}^{l_j}
   R\big(R(z_1)\vec{x}_{\alpha_1}R(z_2)\ldots R(z_j)\vec{x}_{\alpha_j}|_{x_{jt}\to y_{jt}}\ldots R(z_s)\vec{x}_{\alpha_s}x_\beta^{k+1}R(z_{s+1})\big) \bigg).
\end{multline}

In~\eqref{almostRB}--\eqref{LongRel}, we have 
\begin{gather*}
s\geq1,\quad k\geq0,\quad
a,b\in RS(X_\Lambda\cup Y_\Lambda),\quad
z_2,\ldots,z_s\in RS(X_\Lambda\cup Y_\Lambda)\setminus (X_\Lambda\cup Y_\Lambda),\\
\vec{x}_{\alpha_1},\ldots,\vec{x}_{\alpha_{s-1}}\in S(X_\Lambda),\quad
l_j = \deg \vec{x}_{\alpha_j},
\end{gather*}
$\vec{x}_{\alpha_s}\in S(X_\Lambda)$ in~\eqref{XRel} and
$\vec{x}_{\alpha_s}\in S(X_\Lambda)\cup \emptyset$ in~\eqref{LongRel}.
Moreover, $x_\beta$ is greater than all letters from 
$\vec{x}_{\alpha_s} = x_{s1}x_{s2}\ldots x_{s l_s}$.
By $\vec{x}_{\alpha_j}|_{x_{jt}\to y_{jt}}$ we mean the word $\vec{x}_{\alpha_j}$
in which the $t$-th letter $x_{jt}$ is replaced by $y_{jt}$, i.e.,
such $y_{jt}\in Y_\Lambda$ that $R(y_{jt}) = x_{jt}$.

Finally, by $R(z_1)$, we denote either that 
$z_1\in RS(X_\Lambda\cup Y_\Lambda)\setminus (X_\Lambda\cup Y_\Lambda)$ 
or that $R(z_1)$ is absent, i.e., $R(z_1) = \emptyset$.
The same is true for $R(z_{s+1})$.
In particular, for $s=1$, $R(z_1) = R(z_2) = \emptyset$,
the relation \eqref{XRel} equals $R(\vec{x}_{\alpha_1})$.
The values $s=1$, $k=0$, $R(z_1) = R(z_2) = \vec{x}_{\alpha_1} = \emptyset$
transfrom~\eqref{LongRel} to the relation $R(y_\beta) - x_\beta$. 

{\bf Remark}. 
In the relations~\eqref{XRel} and~\eqref{LongRel}, we are using 
associative words $\vec{x}_\alpha\in S(X)$ instead of ordered
polynomials from $\Bbbk[X]$, otherwise we will have to reduce 
the products of such polynomials from $\Bbbk[X]$ to the ordered ones
in all possible compositions from~$S$. 

{\bf Lemma 3.4}.
Given a Lie algebra~$H$, the equality
\begin{equation}\label{Lemma}
(l+1)yx^l
 = \sum\limits_{i=2}^{l+1}(-1)^i \binom{l+1}{i}[y,x^{(i-1)}]x^{l+1-i}
 + \big(yx^l + xyx^{l-1} + \ldots + x^ly)
\end{equation}
holds in the universal enveloping algebra $U(H)$ for any $x,y\in H$ and $l\geq0$.

{\sc Proof}. 
Let us transfer~$y$ to the most left position in the every summand of the sum 
$yx^l + xyx^{l-1} + \ldots + x^ly$.
Thus, we will get by~\eqref{UnivRel} the summands $[y,x^{(i-1)}]x^{l+1-i}$ for $i=2,\ldots,l+1$. 
The coefficient in the sum by the summand 
$(-1)^{i+1} [y,x^{(i-1)}]x^{l+1-i}$ for the fixed~$i$ comes from the sum 
$$
(l+2-i)x^{i-2}[x,y]x^{l-i}+(l+1-i)x^{i-1}[x,y]x^{l-i-1}+\ldots+x^{l-1}[x,y]
$$
and thus equals 
\begin{multline*}
(l+2-i) + (l+1-i)\binom{i-1}{i-2} + (l-i)\binom{i}{i-2} + \ldots + (l+2-i-(l+1-i))\binom{l-1}{i-2} \\
 = (l+1)\left(\binom{i-2}{i-2} + \binom{i-1}{i-2} + \ldots + \binom{l-1}{i-2}\right) \\ \allowdisplaybreaks
 - (i-1)\left( 1 + i + \frac{i(i+1)}{2} + \ldots + \frac{i(i+1)(i+2)\ldots l}{(l-i+1)!}\right) \\
 = (l+1)\binom{l}{i-1} - (i-1)\binom{l+1}{i} = \binom{l+1}{i},
\end{multline*}
and we have proved the formula~\eqref{Lemma}.
\hfill $\Square$

{\bf Theorem 3.5}.
The set $S$ is a a Gr\"{o}bner---Shirshov basis in 
$R\As\langle X_\Lambda\cup Y_\Lambda\rangle$.

{\sc Proof}. 
All compositions between two elements from \eqref{UnivRel} are trivial, 
as it is the method to construct the universal enveloping associative algebra for a given Lie algebra.
Also, compositions of intersection between~\eqref{almostRB} and~\eqref{almostRB} are trivial. 
Thus, all compositions of intersection which are not 
at the same time compositions of inclusion are trivial.

All compositions of inclusion between \eqref{UnivRel} and \eqref{XRel}
are trivial since $[x_\alpha,x_\beta]\in\Span\{X_\Lambda\}$ for all $\alpha,\beta$
(as $\Span\{X_\Lambda\}$ coincides with the image of RB-operator on $\hat{L}$).

Let us compute a composition of inclusion between \eqref{UnivRel} and \eqref{LongRel}.
Let
$$
w = R\big(R(z_1)\vec{x}_{\alpha_1}R(z_2)\ldots R(z_s)\vec{x}_{\alpha_s} y_\beta x_\beta^k R(z_{s+1})\big).
$$
We apply the relation~\eqref{UnivRel}: $x_\gamma x_\delta-x_\delta x_\gamma-[x_\gamma,x_\delta]$ to the subword
$$
\vec{x}_{\alpha_j} = \vec{x}_{\alpha_j}'x_\gamma x_\delta \vec{x}_{\alpha_j}'',\quad x_\gamma>x_\delta,\ 1\leq j\leq s.
$$ 

The triviality of the corresponding composition of inclusion for $j<s$ follows from the equality
\begin{equation}\label{UnivRel-LengRelGen}
R\big(R(z_1)\ldots R(z_j)\vec{x}_{\alpha_j}'
([x_\gamma,y_\delta]+[y_\gamma,x_\delta]-[x_\gamma,x_\delta]_y)\vec{x}_{\alpha_j}''\ldots R(z_s)\vec{x}_{\alpha_s}x_\beta^{k+1}R(z_{s+1})\big) = 0,
\end{equation}
where $[x_\gamma,x_\delta]_y := \sum\limits k_\mu y_\mu$ for 
$[x_\gamma,x_\delta] = \sum\limits k_\mu x_\mu\in\Span\{X_\Lambda\}$.

We have
$$
[x_\gamma,y_\delta]+[y_\gamma,x_\delta]-[x_\gamma,x_\delta]_y
 = [R(y_\gamma),y_\delta]+[y_\gamma,R(y_\delta)]-[R(y_\gamma),R(y_\delta)]_y\in\Span\{X_\Lambda\},
$$
since $R([R(y_\gamma),y_\delta]+[y_\gamma,R(y_\delta)])-[R(y_\gamma),R(y_\delta)] = 0$
in the initial Lie RB-algebra. 
Thus, the expression~\eqref{UnivRel-LengRelGen} is equivalent to zero by~\eqref{XRel}.

Now, study the case $j = s$. 
The triviality of the corresponding composition of inclusion will follow from the equivalence
\begin{multline}\label{UnivRel-LengRelOne}
A(z_1,\vec{x}_{\alpha_1},z_2,\ldots,z_s,\vec{x}_{\alpha_s},z_{s+1}) 
 = R(z_1)\vec{x}_{\alpha_1}R(z_2)\ldots R(z_s)\vec{x}_{\alpha_s}R(z_{s+1}) \\
 - \sum\limits_{j=1}^{s+1}R\big(R(z_1)\vec{x}_{\alpha_1}R(z_2)\ldots R(z_{j-1})\vec{x}_{\alpha_{j-1}}z_j\vec{x}_{\alpha_j}R(z_{j+1})\ldots R(z_s)\vec{x}_{\alpha_s}R(z_{s+1})\big) \\
 - \sum\limits_{j<s,t;\ j=s,t<k_1}
   R\big(R(z_1)\vec{x}_{\alpha_1}R(z_2)\ldots R(z_j)\vec{x}_{\alpha_j}|_{x_{jt}\to y_{jt}}\ldots R(z_s)\vec{x}_{\alpha_s}R(z_{s+1})\big) \bigg) \\
 - R\bigg(R(z_1)\ldots R(z_s)\bigg(\sum\limits_{k_1<t\not\in K}\vec{x}_{\alpha_s}|_{x_{jt}\to y_{jt}}
 + \sum\limits_{t\in K} \big(\vec{x}_{\alpha_s}|_{x_{jt}\to y_{jt}}-\vec{x}_{\alpha_s,0}y_\beta x_\beta^p\big)\bigg) R(z_{s+1})\bigg)  \\
 - (p+1)R\big(R(z_1)\vec{x}_{\alpha_1}R(z_2)\ldots R(z_s)\vec{x}_{\alpha_s,0}y_\beta x_\beta^p R(z_{s+1})\big) \equiv 0 \mod (S,w),
\end{multline}
where $w$ is greater than all terms involved in $A(\ldots)$,
the word $\vec{x}_{\alpha_s} = x_{s1}x_{s2}\ldots x_{s l_s}$ has the biggest letter $x_\beta$
on the positions 
$$
K = \{k_1,k_2,\ldots,k_p\mid k_1<k_2<\ldots<k_p\}\subset\{1,2,\ldots,l_s\}
$$
and the word $\vec{x}_{\alpha_s,0}$ is obtained from $\vec{x}_{\alpha_s}$ by arising all letters $x_\beta$
with preserving order of all remaining letters.

We will proceed on by induction on $l_s = |\vec{x}_{\alpha_s}|$. 
For $l_s = 1$, we are done by~\eqref{LongRel}.

Consider the equality
\begin{equation}\label{GeneralXConvert}
\vec{x}_{\alpha_s}
 = \vec{x}_{\alpha_s,0}x_\beta^p + \sum\limits_{q=1}^p x_{s1}x_{s2}\ldots x_{sk_q-1}\ldots [x_{sk_q},w_q]x_\beta^{p-q},
\end{equation}
where $w_q$ is obtained from the word $x_{sk_q+1}\ldots x_{s l_s}$ by arising all $p-q$ letters $x_\beta$.
For $w_q = w_{q1}w_{q2}\ldots w_{q,l_s-k_q}$, the bracket $[x_{sk_q},w_q]$ in~\eqref{GeneralXConvert} means
$$
[x_{sk_q},w_q] 
 = \sum\limits_{i=1}^{l_s-k_q} w_{q1}\ldots w_{q,i-1}[x_{sk_q},w_{qi}]w_{q,i+1}\ldots w_{q,l_s-k_q}.
$$

By~\eqref{LongRel} and Lemma, we deduce that 
$$
A(z_1,\vec{x}_{\alpha_1},z_2,\ldots,z_s,\vec{x}_{\alpha_s},z_{s+1}) 
 \equiv \sum\limits_{q=1}^p \sum\limits_{i=1}^{l_s-k_q}A(z_1,\vec{x}_{\alpha_1},z_2,\ldots,z_s,\vec{x}_{\alpha_s}(q,i),z_{s+1})\!\! \mod (S,w), 
$$
where 
$\vec{x}_{\alpha_s}(q,i) = x_{s1}\ldots x_{sk_q-1}\ldots w_{q1}\ldots w_{q,i-1}[x_{sk_q},w_{qi}]w_{q,i+1}\ldots w_{q,l_s-k_q}x_\beta^{p-q}$.
The equivalence $A(z_1,\ldots,z_s,\vec{x}_{\alpha_s}(q,i),z_{s+1})\equiv 0$ modulo $(S,w)$ follows by the inductive hypothesis.

Consider a composition of inclusion between \eqref{almostRB} and \eqref{XRel}.
Let $w = R(R(z_1)\vec{x}_{\alpha_1}R(z_2)\ldots$
$R(z_s)\vec{x}_{\alpha_s}R(z_{s+1}))R(b)$
satisfy the above written conditions for~\eqref{XRel}.
At first, we have
$$
R(R(z_1)\vec{x}_{\alpha_1}R(z_2)\ldots R(z_s)\vec{x}_{\alpha_s}R(z_{s+1}))R(b) 
 \mathop{\equiv}\limits^{\eqref{XRel}} 0  \mod (S,w).
$$ 
At second, 
\begin{multline*} 
R(R(z_1)\vec{x}_{\alpha_1}R(z_2)\ldots R(z_s)\vec{x}_{\alpha_s}R(z_{s+1}))R(b) \\
 \mathop{\equiv}\limits^{\eqref{almostRB}}
 R( R(R(z_1)\vec{x}_{\alpha_1}R(z_2)\ldots R(z_s)\vec{x}_{\alpha_s}R(z_{s+1}))b 
 + R(z_1)\vec{x}_{\alpha_1}R(z_2)\ldots R(z_s)\vec{x}_{\alpha_s}R(z_{s+1})R(b)) \\
 \mathop{\equiv}\limits^{\eqref{XRel},\eqref{almostRB}} 
 R( R(z_1)\vec{x}_{\alpha_1}R(z_2)\ldots R(z_s)\vec{x}_{\alpha_s}R(R(z_{s+1})b+z_{s+1}R(b)) ) 
 \mathop{\equiv}\limits^{\eqref{XRel}}  0  \mod (S,w).
\end{multline*}

Compute a composition of inclusion between \eqref{almostRB} and \eqref{LongRel}.
Let $b\in RS(X_\Lambda\cup Y_\Lambda)$ and
$w = R\big(R(z_1)\vec{x}_{\alpha_1} R(z_2)\ldots R(z_s)\vec{x}_{\alpha_s}y_\beta x_\beta^k R(z_{s+1})\big)$
satisfy the above written conditions. Suppose that $R(z_{s+1})\neq\emptyset$, the case $R(z_{s+1}) = \emptyset$ is even easier.
On the one hand, we have
\begin{multline}\label{AlmRB-LongRel:Right}
R\big(R(z_1)\vec{x}_{\alpha_1}R(z_2) \ldots R(z_s)\vec{x}_{\alpha_s}y_\beta x_\beta^k R(z_{s+1})\big)R(b)  \\
 \allowdisplaybreaks
 \mathop{\equiv}\limits^{\eqref{almostRB}} 
R( R\big(R(z_1)\vec{x}_{\alpha_1}R(z_2)\ldots R(z_s)\vec{x}_{\alpha_s}y_\beta x_\beta^k R(z_{s+1})\big)b  ) \\
 + R( R(z_1)\vec{x}_{\alpha_1} R(z_2)\ldots R(z_s)\vec{x}_{\alpha_s}y_\beta x_\beta^k R(z_{s+1}) R(b) ) \\
  \mathop{\equiv}\limits^{\eqref{LongRel},\eqref{almostRB}}
 \frac{1}{k+1} \big( R(R(z_1)\vec{x}_{\alpha_1} R(z_2)\ldots R(z_s)\vec{x}_{\alpha_s}x_\beta^{k+1} R(z_{s+1})b) + R(R(\Sigma) b) \big) \\
 + R( R(z_1)\vec{x}_{\alpha_1} R(z_2)\ldots R(z_s)\vec{x}_{\alpha_s}y_\beta x_\beta^k R(R(z_{s+1})b+z_{s+1}R(b)) ) \mod (S,w),
\end{multline}
where
\begin{multline}\label{Sigma}
\Sigma 
 =  \sum\limits_{i=2}^{k+1}(-1)^i\binom{k+1}{i}R(z_1)\vec{x}_{\alpha_1}R(z_2)\ldots R(z_s)\vec{x}_{\alpha_s}\big[y_\beta,x_\beta^{(i-1)}\big]x_\beta^{k+1-i} R(z_{s+1}) \\
 \allowdisplaybreaks
 - \sum\limits_{j=1}^{s+1}R(z_1)\vec{x}_{\alpha_1} R(z_2)\ldots R(z_{j-1})\vec{x}_{\alpha_{j-1}}z_j\vec{x}_{\alpha_j}R(z_{j+1})\ldots R(z_s)\vec{x}_{\alpha_s}x_\beta^{k+1} R(z_{s+1}) \\
 - \sum\limits_{j=1}^{s}\sum\limits_{t=1}^{l_j}
   R(z_1)\vec{x}_{\alpha_1} R(z_2)\ldots R(z_j)\vec{x}_{\alpha_j}|_{x_{jt}\to y_{jt}}\ldots R(z_s)\vec{x}_{\alpha_s}x_\beta^{k+1}R(z_{s+1}). 
\end{multline}
We proceed on the last summand from~\eqref{AlmRB-LongRel:Right}
\begin{multline*}
R\big( R(z_1)\vec{x}_{\alpha_1} R(z_2)\ldots R(z_s)\vec{x}_{\alpha_s}y_\beta x_\beta^k R(R(z_{s+1})b+z_{s+1}R(b))\big) \\
 \mathop{\equiv}\limits^{\eqref{LongRel}}
 \frac{1}{k+1} \big(R(z_1)\vec{x}_{\alpha_1} R(z_2)\ldots R(z_s)\vec{x}_{\alpha_s}x_\beta^{k+1} R(R(z_{s+1})b+z_{s+1}R(b)) \\
 + R(\Sigma'[R(z_{s+1})\to R(R(z_{s+1})b+z_{s+1}R(b))])  \\
 - R\big(R(z_1)\vec{x}_{\alpha_1} R(z_2)\ldots R(z_s)\vec{x}_{\alpha_s}x_\beta^{k+1}(R(z_{s+1})b+z_{s+1}R(b)) \big) \big)  \mod (S,w)
\end{multline*}
for $\Sigma' = \Sigma + R(z_1)\vec{x}_{\alpha_1} R(z_2)\ldots R(z_s)\vec{x}_{\alpha_s}x_\beta^{k+1}z_{s+1}$.
The sum $\Sigma'[R(z_{s+1})\to R(u)]$ is defined as follows. 
We need to replace the $R$-letter $R(z_{s+1})$ staying at the end by the $R$-letter $R(u)$ in each summand of $\Sigma'$. 

Inserting the obtained expression into~\eqref{AlmRB-LongRel:Right}, we get
\begin{multline}\label{AlmRB-LongRel:RightFin}
 \frac{1}{k+1} \big( 
  R(R(\Sigma) b)  + R(z_1)\vec{x}_{\alpha_1}R(z_2)\ldots R(z_s)\vec{x}_{\alpha_s}x_\beta^{k+1} R(R(z_{s+1})b+z_{s+1}R(b)) \\
 + R(\Sigma'[R(z_{s+1})\to R(R(z_{s+1})b+z_{s+1}R(b))])  \\
 - R\big(R(z_1)\vec{x}_{\alpha_1}R(z_2)\ldots R(z_s)\vec{x}_{\alpha_s}x_\beta^{k+1} z_{s+1}R(b) \big) 
 \big)
\end{multline}

On the other hand, we calculate
\begin{multline*}
R\big(R(z_1)\vec{x}_{\alpha_1}R(z_2)\ldots R(z_s)\vec{x}_{\alpha_s}y_\beta x_\beta^k R(z_{s+1})\big)R(b)  \\
 \mathop{\equiv}\limits^{\eqref{LongRel}} 
\frac{1}{k+1}(R(z_1)\vec{x}_{\alpha_1}R(z_2)\ldots R(z_s)\vec{x}_{\alpha_s}x_\beta^{k+1} R(z_{s+1})+R(\Sigma))R(b) \\
 \mathop{\equiv}\limits^{\eqref{almostRB}} 
 \frac{1}{k+1}\big(
 R(z_1)\vec{x}_{\alpha_1}R(z_2)\ldots R(z_s)\vec{x}_{\alpha_s}x_\beta^{k+1} R(R(z_{s+1})b+z_{s+1}R(b)) \\
 + R(R(\Sigma)b) + R(\Sigma'[R(z_{s+1})\to R(R(z_{s+1})b+z_{s+1}R(b))]) \\
 - R\big(R(z_1)\vec{x}_{\alpha_1} R(z_2)\ldots R(z_s)\vec{x}_{\alpha_s}x_\beta^{k+1} z_{s+1}R(b)\big)
 \big) \mod (S,w). 
\end{multline*}
Thus, the corresponding composition of inclusion is trivial modulo $(S,w)$.

Let us show the triviality of a composition of inclusion between \eqref{LongRel} and \eqref{LongRel}.
Suppose that we have 
\begin{gather*}
w = R\big(R(z_1)\vec{x}_{\alpha_1}R(z_2)\ldots R(z_s)\vec{x}_{\alpha_s}y_\beta x_\beta^k R(z_{s+1})\big),\\
z_m = R(q_1)\vec{x}_{\gamma_1}R(q_2)\ldots R(q_r)\vec{x}_{\gamma_r}y_\delta x_\delta^l R(q_{r+1})
\end{gather*}
for some $1\leq m\leq s+1$. Consider the case $m\leq s$. We will use the denotation $\Sigma$ defined by~\eqref{Sigma} and also
\begin{gather*}\allowdisplaybreaks
\Sigma' := \Sigma + R(z_1)\vec{x}_{\alpha_1}R(z_2)\ldots \vec{x}_{\alpha_{m-1}} z_m\vec{x}_{\alpha_m} \ldots R(z_s)\vec{x}_{\alpha_s}x_\beta^{k+1} R(z_{s+1}), \\
z_m' = R(q_1)\vec{x}_{\gamma_1}R(q_2)\ldots R(q_r)\vec{x}_{\gamma_r}x_\delta^{l+1}R(q_{r+1}), \\
A = \sum\limits_{j=1}^{r+1}R(q_1)\vec{x}_{\gamma_1}R(q_2)\ldots R(q_{j-1})\vec{x}_{\gamma_{j-1}}q_j\vec{x}_{\gamma_j}R(q_{j+1})\ldots R(q_r)\vec{x}_{\gamma_r}x_\delta^{l+1} R(q_{r+1}), \\
B = \sum\limits_{j=1}^r\sum\limits_{t=1}^{m_j}R(q_1)\vec{x}_{\gamma_1} R(q_2)\ldots R(q_j)\vec{x}_{\gamma_j}|_{x_{jt}\to y_{jt}}\ldots R(q_r)\vec{x}_{\gamma_r}x_\delta^{l+1} R(q_{r+1}), 
\quad m_j = \deg \vec{x}_{\gamma_j}, \\ 
\widetilde{\Delta} 
 = -A - B - R(q_1)\vec{x}_{\gamma_1}R(q_2)\ldots R(q_r)\vec{x}_{\gamma_r}\big(y_\delta x_\delta^l+x_\delta y_\delta x_\delta^{l-1}+\ldots +x_\delta^l y_\delta)R(q_{r+1}), \\
\Delta
 =  -A - B + \sum\limits_{i=2}^{l+1}(-1)^i\binom{l+1}{i}R(q_1)\vec{x}_{\gamma_1}R(q_2)\ldots R(q_r)\vec{x}_{\gamma_r} 
 \big[y_\delta,x_\delta^{(i-1)}\big]x_\delta^{l+1-i} R(q_{r+1}). 
\end{gather*}

On the one hand, modulo $(S,w)$ we get
\begin{multline} \label{LongRel-LongRelOut}
R\big(R(z_1)\vec{x}_{\alpha_1}R(z_2)\ldots R(z_s)\vec{x}_{\alpha_s}y_\beta x_\beta^k R(z_{s+1})\big) \\
 \mathop{\equiv}\limits^{\eqref{LongRel},\ \mbox{out}}  
 \frac{1}{k+1}\big(R(z_1)\vec{x}_{\alpha_1}R(z_2)\ldots \vec{x}_{\alpha_{m-1}} R(z_m)\vec{x}_{\alpha_m} \ldots R(z_s)\vec{x}_{\alpha_s}x_\beta^{k+1}R(z_{s+1}) + R(\Sigma)\big) \\
 \mathop{\equiv}\limits^{\eqref{LongRel},\ \mbox{in}}    
 \frac{1}{(k+1)(l+1)}\big(R(z_1)\vec{x}_{\alpha_1}R(z_2)\ldots \vec{x}_{\alpha_{m-1}}z_m'\vec{x}_{\alpha_m}\ldots R(z_s)\vec{x}_{\alpha_s}x_\beta^{k+1} R(z_{s+1}) \\ 
 + R(z_1)\vec{x}_{\alpha_1}R(z_2)\ldots \vec{x}_{\alpha_{m-1}}R(\Delta)\vec{x}_{\alpha_m}\ldots R(z_s)\vec{x}_{\alpha_s}x_\beta^{k+1} R(z_{s+1})
 + R(\Sigma' |_{R(z_m) \to z_m'}) \\
 + R(\Sigma' |_{R(z_m)\to R(\Delta)})
 - (l+1)R\big(R(z_1)\ldots \vec{x}_{\alpha_{m-1}}z_m\vec{x}_{\alpha_m}\ldots R(z_s)\vec{x}_{\alpha_s}x_\beta^{k+1} R(z_{s+1})\big)\big).
\end{multline}

On the other hand, 
\begin{multline}\label{LongRel-LongRelIn} 
R\big(R(z_1)\vec{x}_{\alpha_1}R(z_2)\vec{x}_{\alpha_2}\ldots R(z_s)\vec{x}_{\alpha_s}y_\beta x_\beta^k R(z_{s+1})\big) \\
 \mathop{\equiv}\limits^{\eqref{LongRel},\ \mbox{in}}  
 \frac{1}{l+1}R\big(R(z_1)\vec{x}_{\alpha_1}R(z_2)\ldots \vec{x}_{\alpha_{m-1}} z_m'\vec{x}_{\alpha_m} \ldots R(z_s)\vec{x}_{\alpha_s}y_\beta x_\beta^k R(z_{s+1}) \\
 + R\big(R(z_1)\vec{x}_{\alpha_1}R(z_2)\ldots \vec{x}_{\alpha_{m-1}} R(\Delta)\vec{x}_{\alpha_m} \ldots R(z_s)\vec{x}_{\alpha_s}y_\beta x_\beta^k R(z_{s+1})\big) \big) \\
 \mathop{\equiv}\limits^{\eqref{LongRel},\ \mbox{out}}    
 \frac{1}{(k+1)(l+1)}\big(R(z_1)\vec{x}_{\alpha_1}R(z_2)\ldots \vec{x}_{\alpha_{m-1}}z_m'\vec{x}_{\alpha_m}\ldots R(z_s)\vec{x}_{\alpha_s}x_\beta^{k+1} R(z_{s+1}) \\ 
  + R(\Sigma' |_{R(z_m) \to z_m'}) 
  - R\big(R(z_1)\vec{x}_{\alpha_1}R(z_2)\ldots \vec{x}_{\alpha_{m-1}}\widetilde{\Delta}\vec{x}_{\alpha_m}\ldots R(z_s)\vec{x}_{\alpha_s}x_\beta^{k+1} R(z_{s+1})\big) \\
  + R(z_1)\ldots \vec{x}_{\alpha_{m-1}}R(\Delta)\vec{x}_{\alpha_m}\ldots R(z_s)\vec{x}_{\alpha_s}x_\beta^{k+1} R(z_{s+1}) + R(\Sigma |_{z_m\to \Delta})\big)
   \mod (S,w).
\end{multline}

Subtracting~\eqref{LongRel-LongRelIn} from~\eqref{LongRel-LongRelOut}, we have up to the factor $(k+1)(l+1)$
$$
R\big(R(z_1)\vec{x}_{\alpha_1}R(z_2)\ldots \vec{x}_{\alpha_{m-1}}
C \vec{x}_{\alpha_m}\ldots R(z_s)\vec{x}_{\alpha_s}x_\beta^{k+1}R(z_{s+1})\big),
$$
where
\begin{multline*}
C = \Delta - \widetilde{\Delta} -  (l+1)R(q_1)\vec{x}_{\gamma_1}R(q_2)\ldots R(q_r)\vec{x}_{\gamma_r}y_\delta x_\delta^l R(q_{r+1}) \\
  = R(q_1)\vec{x}_{\gamma_1}R(q_2)\vec{x}_{\gamma_2}\ldots R(q_r)\vec{x}_{\gamma_r}DR(q_{r+1})
\end{multline*}
for 
$$
D = \sum\limits_{i=2}^{l+1}(-1)^i \binom{l+1}{i}[y_\delta,x_\delta^{(i-1)}]x_\delta^{l+1-i}
 + \big(y_\delta x_\delta^l+x_\delta y_\delta x_\delta^{l-1}+\ldots + x_\delta^l y_\delta) - (l+1)y_\delta x_\delta^l.
$$
Equality of $D$ to zero follows from Lemma 3.4.

The proof in the case~$m = s+1$ is quite similiar. We only need one additional thing: we should apply~\eqref{UnivRel-LengRelOne} 
instead of~\eqref{LongRel} for the term  
$R\big(R(z_1)\ldots R(z_s)\vec{x}_{\alpha_s}y_\beta x_\beta^k z_m'\big)$.
Such application is correct, since all terms involved in~\eqref{UnivRel-LengRelOne} 
have less $R$-degree than $w$.

It is easy to verify that all other compositions of inclusion are trivial.
\hfill $\Square$

{\bf Corollary 3.6}.
The quotient $A$ of $R\As\langle X_\Lambda\cup Y_\Lambda\rangle$ by $Id(S)$
is the universal enveloping associative RB-algebra for the Lie algebra 
$\hat{L}$ with the RB-operator $R$.
Moreover, $\hat{L}$ injectively embeds into $A^{(-)}$.

{\sc Proof}. 
By~\eqref{almostRB}, $A$ is an associative RB-algebra.
By \eqref{UnivRel}--\eqref{XRel}, 
we have that $A$ is enveloping of $\hat{L}$
for both: the Lie bracket $[,]$ and the action of $R$.
Thus, $A$ is an associative enveloping of $\hat{L}$.

Let us prove that $A$ is the universal enveloping one.
At first, $A$ is generated by~$\hat{L}$. At second,
all elements from $S$ are identities in the universal enveloping
associative RB-algebra $U_{RB}(\hat{L})$. Indeed,
\eqref{UnivRel} are enveloping conditions for the product,
\eqref{almostRB} is the RB-identity,
\eqref{XRel} and \eqref{LongRel} for $s=1$, $k=0$, $R(z_1) = R(z_2) = \vec{x}_{\alpha_1} = \emptyset$
give the relations~\eqref{R-Cond}, the enveloping conditions for the action of $R$ on $L'$.
By~\eqref{R-Cond} and~\eqref{almostRB},
we deduce that $R(\vec{x}_\alpha) = 0$ for $\vec{x}_\alpha\in S(X)$.

Let us prove the equality
$$
w = R(R(z_1)\vec{x}_{\alpha_1}R(z_2)\vec{x}_{\alpha_2}\ldots R(z_s)\vec{x}_{\alpha_s}R(z_{s+1})) = 0.
$$
The relation 
\begin{multline}\label{LongRB}
R(a_1)R(a_2)\ldots R(a_k) 
 = R(a_1R(a_2)\ldots R(a_k)+R(a_1)a_2R(a_3)\ldots R(a_k) \\
 +\ldots+R(a_1)R(a_2)\ldots R(a_{k-1})a_k)
\end{multline}
holds in an associative RB-algebra as the direct consequence of~\eqref{RB}.

By~\eqref{LongRB}, we may rewrite $w$ as $R(R(a)\vec{x}_{\alpha_s}R(z_{s+1}))$
for some $a$, since $x = R(y)$ for each $x\in L$. Now we again by~\eqref{LongRB} get 
\begin{multline*}
w = R(R(a)\vec{x}_{\alpha_s}R(z_{s+1})) \\
  = R(a)R(\vec{x}_{\alpha_s})R(z_{s+1}) 
  - R(aR(\vec{x}_{\alpha_s})R(z_{s+1}))
  - R(R(a)R(\vec{x}_{\alpha_s})z_{s+1})
  = 0
\end{multline*}
proving that~\eqref{XRel} holds in~$U_{RB}(\hat{L})$.

The relation~\eqref{LongRel} follows from~\eqref{LongRB} and the next one fullfiled in~$U_{RB}(\hat{L})$:
\begin{equation}\label{YXRel}
R(yx^l)
 = \frac{1}{l+1}\left(   
 x^{l+1} + \sum\limits_{i=2}^{l+1}(-1)^i \binom{l+1}{i}R([y,x^{(i-1)}]x^{l+1-i}) \right).
\end{equation}

Let us state~\eqref{YXRel}. By~\eqref{LongRB}, we have
\begin{equation}\label{ProofOfYXRel}
x^{l+1} = R(y)R(y)\ldots R(y) 
 = R(yx^l + xyx^{l-1} + x^2 yx^{l-2} + \ldots + x^l y).
\end{equation}
It remains to apply Lemma~3.4.
\hfill $\Square$

Finally, by Theorems 3.3 and 3.5 we get the injectivity of
embedding $\hat{L}$ into $A^{(-)}$.

{\bf Corollary 3.7}.
Any pre-Lie algebra injectively embeds into its 
universal enveloping preassociative algebra.

{\sc Proof}. 
Let $L$ be a pre-Lie algebra. By Theorem~2.2, 
$L$ can be injectively embedded into $\hat{L}^{(R)}$
with the RB-operator $R$ of weight~0. 
Then, by Corollary~3.6, we embed the Lie RB-algebra $\hat{L}$ 
into its universal enveloping associative algebra $A$ with the RB-operator~$P$. 
Thus, the subalgebra (in prealgebra sense) $T$ in $A^{(P)}$ generated by the set~$L'$ 
is an (injective) enveloping preassociative algebra of initial pre-Lie algebra $L$.
\hfill $\Square$

\section*{Acknowledgements}

This work was supported by the Austrian Science Foundation FWF grant P28079.

\noindent Vsevolod Gubarev \\
University of Vienna \\
Oskar-Morgenstern-Platz 1, 1090, Vienna, Austria \\
e-mail: vsevolod.gubarev@univie.ac.at

\end{document}